\documentclass[titlepage,draft,12pt]{article} 
\usepackage{amsfonts,latexsym,pstricks,pst-plot} 
\usepackage[a4paper]{geometry}

\newcommand{\re}{{\mathbb{R}}}
\newcommand{\n}{{\mathbb{N}}}
\newcommand{\ep}{\varepsilon}

\newtheorem{thm}{Theorem}[section]

\title{On the evolution of subcritical regions for the Perona-Malik
equation}

\author{Marina Ghisi\vspace{1ex}\\ {\normalsize Universit\`a degli
Studi di Pisa} \\{\normalsize Dipartimento di Matematica ``Leonida
Tonelli''}\\
{\normalsize PISA (Italy)}\\
{\normalsize e-mail: \texttt{ghisi@dm.unipi.it}}\and
Massimo Gobbino\vspace{1ex}\\ {\normalsize Universit\`a degli Studi di Pisa} 
\\{\normalsize Dipartimento di Matematica Applicata ``Ulisse Dini''}\\ 
{\normalsize 
 PISA (Italy)}\\  
{\normalsize e-mail: \texttt{m.gobbino@dma.unipi.it}}}

\date{}

\begin{document}
\maketitle
\begin{abstract}
	
	The Perona-Malik equation is a celebrated example of
	forward-backward parabolic equation.  The forward behavior takes
	place in the so-called subcritical region, in which the gradient
	of the solution is smaller than a fixed threshold.  In this paper
	we show that this subcritical region evolves in a different way in
	the following three cases: dimension one, radial solutions in
	dimension greater than one, general solutions in dimension greater
	than one.
	
	In the first case subcritical regions increase, but there is no
	estimate on the expansion rate.  In the second case they expand
	with a positive rate and \emph{always spread over the whole
	domain} after a finite time, depending only on the (outer) radius
	of the domain.  As a by-product, we obtain a non-existence result
	for global-in-time classical radial solutions with large enough
	gradient.  In the third case we show an example where subcritical
	regions do \emph{not} expand.
	
	Our proofs exploit comparison principles for suitable degenerate
	and non-smooth free boundary problems.
	
\vspace{1cm}

\noindent{\bf Mathematics Subject Classification 2000 (MSC2000):}
35K55, 35K65, 35R35.

\vspace{1cm} \noindent{\bf Key words:} Perona-Malik equation,
forward-backward parabolic equation, degenerate parabolic equation,
moving domains, subsolutions and supersolutions.
\end{abstract}
 
\section{Introduction}

In this paper we consider the Perona-Malik equation
\begin{equation}
	u_{t}(x,t)-\mathrm{div}\left(\frac{\nabla u(x,t)}{1+|\nabla
	u(x,t)|^{2}}\right)  =  0
	\quad\quad
	\forall (x,t)\in\Omega\times[0,T),
	\label{eq:PM} 
\end{equation}
where $\Omega\subseteq\re^{n}$ is an open set and $T>0$. This equation
is the formal gradient flow of the functional
$$PM(u):=\frac{1}{2}\int_{\Omega}^{}\log\left(1+|\nabla
u(x)|^{2}\right)dx.$$

The convex-concave behavior of the integrand makes (\ref{eq:PM}) a
forward-backward partial differential equation of parabolic type.  The
forward (or \emph{subcritical}) region is the set of points $(x,t)$
where $|\nabla u(x,t)|<1$, the backward (or \emph{supercritical})
region is the set of points where $|\nabla u(x,t)|>1$.

This equation was introduced by \textsc{P.\ Perona} and \textsc{J.\
Malik}~\cite{PM} in the context of image denoising. They considered
equation (\ref{eq:PM}) in a rectangle $\Omega\subseteq\re^{2}$, with
an initial condition $u(x,0)=u_{0}(x)$ representing the grey level of 
a (noisy) picture, and Neumann boundary conditions. For increasing
values of $t$ the functions $x\to u(x,t)$ are thought as successively
restored versions of $u_{0}(x)$.

The heuristic idea is that small disturbances, corresponding to small
values of the gradient, are smoothed out because of the diffusion
which takes place in the forward regions.  On the other hand, sharp
edges correspond to large values of the gradient and therefore they
are expected to be enhanced by the backward character of the equation
in supercritical regions.  This phenomenology has been actually
observed in numerical experiments, which also reveal an unexpected
stability (see~\cite{E2,E3,gl}).  This discrepancy between the
practical efficacy of (\ref{eq:PM}) and its analytical ill-posedness
has been called ``Perona-Malik paradox'' after \textsc{S.\
Kichenassamy}'s paper~\cite{kich}.

In the last fifteen years the paradox has been investigated in
numerous papers.  Several authors proved well posedness results for
approximations of (\ref{eq:PM}) obtained via space
discretization~\cite{BNPT1,GG} or convolution~\cite{CLMC}, time
delay~\cite{amann}, fractional derivatives~\cite{guidotti}, fourth
order regularization~\cite{BF}, simplified nonlinearities~\cite{BNP}.
The behavior of such approximations as the suitable parameter goes to
0 is a much more challenging problem.  As far as we know, results in
this direction have been obtained only for the semidiscrete scheme in
dimension one.  In this case the authors~\cite{GG} (see also
\cite{BNPT1}) proved that approximated solutions converge to a limit,
and under reasonable assumptions this limit is a classical solution of
(\ref{eq:PM}) inside its subcritical region.

Despite of these partial progresses, a solution of the paradox is
still far away.  We remind that a solution of the paradox is a notion
of weak solution for (\ref{eq:PM}) which exists for large classes of
initial data (for example in $BV(\Omega)$ or $SBV(\Omega)$), is
reasonably stable, and to which reasonable approximations converge.

In this direction, \textsc{K.\ Zhang}~\cite{Z1} (see also
\cite{Z2,Z3}) showed that the class of Lipschitz solutions is far from
being a solution to the paradox.  Indeed he proved that for any
nonconstant smooth initial condition, even if subcritical, the Neumann
boundary value problem admits infinitely many (pathological) Lipschitz
solutions.

Classical solutions (namely solutions which are at least of class
$C^{1}$) have also been investigated in the last decade.  \textsc{B.\
Kawohl} and \textsc{N.\ Kutev}~\cite{kk} observed that global-in-time
classical solutions exist if the initial condition is subcritical,
while in~\cite{kich} it is remarked that local-in-time classical
solutions cannot exist unless the initial condition is very regular in
its supercritical region.  Moreover the authors proved in
\cite{GG-LocSol} that in dimension one there exists a dense set of
initial data for which the Cauchy problem with Neumann boundary
conditions admits a local-in-time classical solution of class
$C^{2,1}$ (namely with two continuous derivatives with respect to
space variables, and one continuous derivative with respect to time).
On the other hand, such solutions \emph{cannot} be global if the
initial condition has a nonempty supercritical region (see \cite{kk}
and \cite{G-PMEntire}).

Quite surprisingly, things are not so drastic in dimension greater
than one. Indeed the authors proved in~\cite{GG-GlobSol} that
\emph{global-in-time} radial solutions of class $C^{2,1}$ \emph{do}
exist for some classes of initial data with nonempty supercritical
region.

For classical solutions one can define the family of open sets
$$I^{-}(t):=\left\{x\in\Omega:\left|\nabla u(x,t)\right|
<1\right\}
\quad\quad
\forall t\in[0,T).$$

This is the family of subcritical regions.  Its behavior as $t$ varies
is the object of this paper.  We point out that this definition is
purely local, in the sense that it does not depend on the boundary
conditions.  We show three situations in which subcritical regions
evolve in a different way.  As far as we know, these different
behaviors had not been explicitly reported in numerical experiments.
The reason is maybe that in an unstable framework it is always
difficult to distinguish what depends on the model itself, and what
depends on the implementation.  In any case, we leave to numerical
analysts and applied mathematicians any comment about the significance
of these results for the model and its practical applications.

\paragraph{\textmd{\emph{The one dimensional case}}}

In dimension one subcritical regions don't shrink, namely
\begin{equation}
	I^{-}(s)\subseteq I^{-}(t)
	\quad\quad
	\mbox{whenever } 0\leq s\leq t<T.
	\label{shrinking}
\end{equation}

This fact had already been proved in \cite{kk} under some structural
assumptions on the initial condition, afterwards removed in
\cite{G-PMEntire}.

In this paper we present an alternative proof (see Theorem~\ref{thm:PM1})
based on a comparison principle, which only requires $C^{1}$
regularity, and works substantially for all equations which are the
gradient flow of a nonconvex functional.  This proof gives us the
opportunity to show, in a simpler setting, the method which is
fundamental in the next case, when the result was not known before.

This result is optimal in the sense that it may happen that the
subcritical region is the same for every $t\in[0,T)$.  An example is
provided by the local-in-time solutions constructed
in~\cite{GG-LocSol}.

\paragraph{\textmd{\emph{The radial case}}}

Let us consider a radial solution of (\ref{eq:PM}) defined in a ball
or an anulus.  Then (\ref{shrinking}) holds true also in this case.
If moreover $I^{-}(0)\neq\emptyset$, then the inclusion is strict
whenever $s<t$, and there exists $T_{0}>0$ such that $I^{-}(t)=\Omega$
for every $t>T_{0}$.  The value of $T_{0}$ depends on the outer radius
of $\Omega$, but is independent on $u$.  In other words, supercritical
and critical regions disappear after a finite time depending only on
$\Omega$.

If the solution survives and remains regular up to $T_{0}$, then it
becomes subcritical and there are no more obstructions to global
existence. This is actually what happens in the classes of
global-in-time radial solutions constructed in~\cite{GG-GlobSol}.

In Theorem~\ref{thm:no-ex} we show that this is not always the case.
If the initial condition has a supercritical region where the gradient
is large enough, and this supercritical region is surrounded by
subcritical regions, then there is no $C^{1}$ classical solution with
$T>T_{0}$, independently on the boundary conditions.  The reason is
that the maximum of $|\nabla u(x,t)|$ in this supercritical region is
a function of time which cannot decrease too fast.  On the other hand,
the supercritical region must disappear after $T_{0}$, and thus this
maximum has not enough time to decrease from its large initial value
up to 1.

\paragraph{\textmd{\emph{The nonradial case}}}

Roughly speaking, in the radial case each interface between the
subcritical and the supercritical region is a circle which evolves with
velocity proportional to the inverse of its radius.  This reminded us
of the mean curvature motion, and in a first moment led us to suspect
that this interface could evolve in a similar way also for nonradial
solutions.  In particular it seemed reasonable that this interface
could evolve in such a way that supercritical regions tend to shrink,
at least where they are convex.

In Theorem~\ref{thm:example} we show that this is not the case.  We
prove indeed that there exists an initial condition $u_{0}$ in two
variables with the following properties.  The supercritical region of
$u_{0}$ is convex in a neighborhood of the origin, and any
local-in-time solution starting from $u_{0}$, independently on the
boundary conditions, has a supercritical region which expands in a
neighborhood of the origin.  In particular any such solution does
\emph{not} satisfy (\ref{shrinking}).

\paragraph{\textmd{\emph{Evolution of supercritical regions}}}

If subcritical regions expand, then supercritical regions shrink.
This is true, but not obvious.  Indeed it could happen that both the
subcritical and the supercritical region expand at the expense of the
critical region where $|\nabla u(x,t)|=1$.  Nevertheless one can prove
that this is not the case.  It is enough to apply the techniques of
this paper to the equation obtained by reversing the time.  In this
way the role of subcritical and supercritical regions is just
exchanged.

\paragraph{\emph{\textmd{Connection with free boundary problems}}}
 
The evolution of subcritical regions is itself a free boundary
problem.  Unfortunately it involves a forward-backward equation.
After some variable changes we reduce ourselves to more standard
situations. Roughly speaking, in the new variable $v$ we end up with
equations such as
$$v_{t}=\sqrt{v}\,v_{xx},
\hspace{3em}
v_{t}=\sqrt{v}\left\{v_{xx}+A+\mbox{ lower order terms}\right\},$$
where $A>0$. These equations are satisfied where $v>0$, and we are
interested in the evolution of the region where $v$ is positive.

The good news are that these equations are forward parabolic. The bad 
news are that they are \emph{degenerate} and they involve a nonlinear 
term which is \emph{not Lipschitz} continuous. This complicates things
when using comparison principles.

For the first equation we show (see Theorem~\ref{thm:FBP1}) that the
region where $v$ is positive does not shrink with time.  This is
enough to threat the Perona-Malik equation in dimension one.  For the
second equation we show (see Theorem~\ref{thm:FBP2}) that the region
where $v$ is positive expands with a positive rate depending on $A$.
This implies our conclusions for the radial Perona-Malik
equation.\bigskip

This paper is organized as follows.  In section~\ref{sec:statements}
we state our results for the Perona-Malik equation and the related
free boundary problems. In section~\ref{sec:proofs} we prove these
results.  

\setcounter{equation}{0}
\section{Statements}\label{sec:statements}

Throughout this paper we assume that $\varphi\in C^{\infty}(\re)$ is
an even function, hence in particular $\varphi'(0)=0$. We also assume 
that
\begin{equation}
	\varphi''(\sigma)>0
	\hspace{2em}
	\forall\sigma\in[0,1),
	\label{hp:phi+}
\end{equation}
\begin{equation}
	\varphi''(1)=0,
	\label{hp:phi=}
\end{equation}
\begin{equation}
	\varphi''(\sigma)<0
	\hspace{2em}
	\forall\sigma> 1.
	\label{hp:phi-}
\end{equation}

These assumptions imply that $\varphi'(1)>0$ and $\varphi'''(1)\leq 0$.
In some statement we also need the stronger assumption
\begin{equation}
	\varphi'''(1)<0.
	\label{hp:phi'''}
\end{equation}

These assumptions are consistent with the concrete case
$\varphi(\sigma)=2^{-1}\log(1+\sigma^{2})$ of the Perona-Malik
equation. We refer to Figure~\ref{fig:h} for the typical behavior of
$\varphi'(\sigma)$.

We consider the following equation
\begin{equation}
	u_{t}=\mathrm{div}\left(\varphi'\left(|\nabla u|\right)
	\frac{\nabla u}{|\nabla u|}\right),
	\label{eq:PMn}
\end{equation}
which is the natural generalization of (\ref{eq:PM}).  We believe and
we hope that this generality simplifies the presentation, and shows
more clearly which properties of the nonlinearity are essential in
each step.  For the sake of generality one could also weaken the
regularity assumptions on $\varphi$ (we never consider more than 3
derivatives), replace the threshold $\sigma=1$ in (\ref{hp:phi=}) with
any $\sigma_{1}$, and weaken (\ref{hp:phi+}) and (\ref{hp:phi-}) by
asking only that $\varphi''$ is positive in a left-hand neighborhood
of $\sigma_{1}$ and negative a right-hand neighborhood of
$\sigma_{1}$.

\subsection{Main results}

Let us state our results on the evolution of subcritical regions.  The
first result concerns the one dimensional case, where (\ref{eq:PMn})
reduces to
\begin{equation}
	u_{t}=\left(\varphi'(u_{x})\right)_{x}
	\label{eq:PM1}
\end{equation}

This form of the equation is suitable for $C^{1}$ solutions, because
it involves only first order derivatives. When the solution is of
class $C^{2,1}$, equation (\ref{eq:PM1}) can of course be rewritten as
$u_{t}=\varphi''(u_{x})u_{xx}$.

\begin{thm}\label{thm:PM1}
	Let $\varphi\in C^{\infty}(\re)$ be a function satisfying
	(\ref{hp:phi+}) through (\ref{hp:phi-}). 
	Let $x_{1}\leq x_{3}< x_{4}\leq x_{2}$ and $T>0$ be real numbers. 
	
	Let $u\in C^{1}\left((x_{1},x_{2})\times[0,T)\right)$ be a
	function satisfying (\ref{eq:PM1}) in $(x_{1},x_{2})\times[0,T)$,
	and
	\begin{equation}
		|u_{x}(x,0)|<1
		\quad\quad
		\forall x\in(x_{3},x_{4}).
		\label{hp:PM1}
	\end{equation}
	
	Then $|u_{x}(x,t)|<1$ for every
	$(x,t)\in(x_{3},x_{4})\times[0,T)$.
\end{thm}

We point out that in Theorem~\ref{thm:PM1} we don't need assumption
(\ref{hp:phi'''}), which was used in all previous results of the same
type.

Our second result concerns the radial case. Let $r:=|x|$ be the radial 
variable, and let $u(r,t)$ be a radial solution. In this case
(\ref{eq:PM}) becomes
\begin{equation}
	u_{t}=\left(\varphi'(u_{r})\right)_{r}+
	(n-1)\frac{\varphi'(u_{r})}{r},
	\label{eq:PMR}
\end{equation}
where $n$ is the space dimension. From now on we assume for simplicity
that $n=2$. The general case is completely analogous. 

\begin{thm}\label{thm:PMR}
	Let $\varphi\in C^{\infty}(\re)$ be a function satisfying
	(\ref{hp:phi+}) through (\ref{hp:phi'''}).  Let
	$0<r_{1}\leq r_{3}< r_{4}\leq r_{2}$ and $T>0$ be real numbers.  Let
	$u\in C^{1}\left((r_{1},r_{2})\times[0,T)\right)$ be a function
	satisfying (\ref{eq:PMR}) (with $n=2$) in
	$(r_{1},r_{2})\times[0,T)$, and 
	$$|u_{r}(r,0)|<1 \quad\quad
	\forall r\in(r_{3},r_{4}).$$
	
	Let $k_{0}:=r_{2}^{-1}\sqrt{2\varphi'(1)|\varphi'''(1)|}$, and let
	us set
	\begin{equation}
		\mathcal{D}:=\{(r,t)\in(r_{1},r_{2})\times[0,T):
		r_{3}-k_{0}t<r<r_{4}+k_{0}t\}.
		\label{defn:D}
	\end{equation}
	
	Then $|u_{r}(r,t)|<1$ for every $(r,t)\in\mathcal{D}$.
\end{thm}

In other words, this result says that in the radial case the
subcritical region expands with a rate which is bounded from below by
a positive constant $k_{0}$.  Figure~\ref{fig:D-star} shows the shape
of the set $\mathcal{D}$.  The slope of the slanted lines depends on
$k_{0}$.  It is clear that when $t>(r_{2}-r_{1})/k_{0}$ every
nonempty initial subcritical region $(r_{3},r_{4})$ has invaded the
whole interval $(r_{1},r_{2})$.

The third result concerns the nonexistence of global-in-time $C^{1}$
radial solutions if the gradient of the initial condition is too
large.  We point out that this result is independent on the boundary
conditions.

\begin{thm}\label{thm:no-ex}
	Let $\varphi\in C^{\infty}(\re)$ be a function satisfying
	(\ref{hp:phi+}) through (\ref{hp:phi'''}).  Let
	$0\leq r_{1}<r_{3}<r_{4}<r_{5}<r_{2}$ and $T>0$ be real numbers,
	and let $u\in C^{1}\left((r_{1},r_{2})\times[0,T)\right)$ be a
	solution of (\ref{eq:PMR}) (with $n=2$) such that
	\begin{equation}
		|u_{r}(r_{3},0)|<1,
		\hspace{1.5em}
		|u_{r}(r_{5},0)|<1,
		\hspace{1.5em}
		u_{r}(r_{4},0)>
		1+\frac{r_{2}(r_{2}-r_{1})}{r_{1}^{2}}\sqrt{\frac{\varphi'(1)}{
		2|\varphi'''(1)|}}.
		\label{hp:i-data}
	\end{equation}
	
	Then we have that
	$$T\leq\frac{r_{2}(r_{2}-r_{1})}{\sqrt{2\varphi'(1)|\varphi'''(1)|}}.$$
\end{thm}

Our last result is a counterexample to the expansion of subcritical
regions for nonradial solutions.  Note that the condition we impose on
the initial condition $u_{0}$ depends only on the Taylor expansion of
$u_{0}$ of order 3 in a neighborhood of the origin.

\begin{thm}\label{thm:example}
	Let $\varphi\in C^{\infty}(\re)$ be a function such that
	$\varphi'(1)>0$ and $\varphi''(1)=0$.  Let
	$\Omega\subseteq\re^{2}$ be any open set such that
	$(0,0)\in\Omega$.  Let $u_{0}:\Omega\to\re$ be any function of
	class $C^{3}$ such that
	\begin{equation}
		u_{0}(x,y)=\frac{\sqrt{2}}{2}x+\frac{\sqrt{2}}{2}y+
		k_{1}x^{2}+k_{2}y^{2}+
		h_{1}x^{3}+h_{2}y^{3}+
		o\left((x^{2}+y^{2})^{3/2}\right)
		\label{hp:taylor}
	\end{equation}
	as $(x,y)\to(0,0)$.
	This clearly implies that $|\nabla u_{0}(0,0)|=1$. Moreover one
	can choose the parameters $k_{1}$, $k_{2}$, $h_{1}$, $h_{2}$ in
	such a way that the following properties hold true.
	\begin{enumerate}
		\renewcommand{\labelenumi}{(\arabic{enumi})}
		\item  There exist $\delta>0$, $a>0$, and a \emph{convex}
		function $g:(-\delta,\delta)\to(-a,a)$ such that
		$$\hspace{-3em}\left\{(x,y)\in(-\delta,\delta)\times(-a,a):
		|\nabla u_{0}(x,y)|>1\right\}=$$
		$$\hspace{3em}=\left\{(x,y)\in(-\delta,\delta)\times(-a,a):
		y>g(x)\right\}.$$
	
		\item  Let $T>0$, and let $u\in C^{2}(\Omega\times[0,T))$ be a
		function satisfying (\ref{eq:PMn}), and the initial condition
		$u(x,y,0)=u_{0}(x,y)$ for every $(x,y)\in\Omega$. Then
		there exists $t_{0}\in(0,T)$ such that
		$$|\nabla u(0,0,t)|>1
		\quad\quad
		\forall t\in(0,t_{0}).$$
	\end{enumerate}
\end{thm}

In other words, at time $t=0$ the origin lies in the interface which
separates the subcritical and the (locally convex) supercritical
region, while for $t\in(0,t_{0})$ the origin has been absorbed by the
supercritical region.

\subsection{Heuristics}\label{sec:heuristic}

In this section we present simple ``proofs'' of Theorem~\ref{thm:PM1} 
and Theorem~\ref{thm:PMR}.

Let us start with Theorem~\ref{thm:PM1}.  Let us assume that $u$ is
smooth enough, and that one component of the interface between the
subcritical and the supercritical region is represented by a smooth
curve $(\alpha(t),t)$.  Just to fix ideas, let us assume that the
subcritical region lies on the left of the interface, namely where
$r<\alpha(t)$.  Taking the time derivative of the relation
$u_{x}(\alpha(t),t)=1$ we obtain that
\begin{eqnarray*}
	0 & = & u_{xx}(\alpha(t),t)\alpha'(t)+u_{xt}(\alpha(t),t)  \\
	 & = & u_{xx}\alpha'(t)+\varphi''(u_{x})u_{xxx}+
	 \varphi'''(u_{x})u_{xx}^{2},
\end{eqnarray*}
where all the partial derivatives of $u$ in the last line are computed
in the point $(\alpha(t),t)$.  Recalling that
$\varphi''(u_{x}(\alpha(t),t))=0$, we have therefore that
$$\alpha'(t)=-\varphi'''(1)u_{xx}(\alpha(t),t).$$

Now we have that $\varphi'''(1)\leq 0$, and $u_{xx}(\alpha(t),t)\geq
0$ because the subcritical region lies on the left of $r=\alpha(t)$.
We have thus proved that $\alpha'(t)\geq 0$, hence the subcritical
region tends to expand.

This ``proof'' is \emph{not} rigorous for several reasons: we assumed
that $u$ is of class $C^{3}$, we assumed that the interface is a
smooth curve, we divided by $u_{xx}$ which could be 0.  Nevertheless
we believe that this simple argument is quite explicative.

Let us consider now Theorem~\ref{thm:PMR}. As before, we assume that
the interface is given by a smooth curve $(\alpha(t),t)$, the
subcritical region being on the left. Taking the time
derivative of the relation $u_{r}(\alpha(t),t)=1$ we obtain that
\begin{eqnarray*}
	0 & = & u_{rr}(\alpha(t),t)\alpha'(t)+u_{rt}(\alpha(t),t)  \\
	 & = & u_{rr}\alpha'(t)+\varphi''(u_{r})u_{rrr}+
	 \varphi'''(u_{r})u_{rr}^{2}+
	 \frac{\varphi''(u_{r})}{\alpha(t)}u_{rr}-
	 \frac{\varphi'(u_{r})}{\alpha^{2}(t)},
\end{eqnarray*}
where all the partial derivatives of $u$ in the last line are computed
in the point $(\alpha(t),t)$.  Recalling that
$\varphi''(u_{r}(\alpha(t),t))=0$, we have therefore that
$$\alpha'(t)=\frac{\varphi'(1)}{\alpha^{2}(t)}\cdot\frac{1}{u_{rr}}
-\varphi'''(1)u_{rr}= \frac{\varphi'(1)}{\alpha^{2}(t)}
\cdot\frac{1}{u_{rr}} +|\varphi'''(1)|u_{rr}.$$

Applying the inequality between the arithmetic and the geometric mean,
we thus obtain that
$$\alpha'(t)\geq
\frac{2\sqrt{\varphi'(1)|\varphi'''(1)|}}{\alpha(t)}\geq
\frac{2\sqrt{\varphi'(1)|\varphi'''(1)|}}{r_{2}}.$$

This ``proves'' that the subcritical region expands with a rate which 
is bounded from below by a positive constant. The value of this
constant is quite similar to the constant $k_{0}$ of
Theorem~\ref{thm:PMR}.

\subsection{Free boundary problems}

Our proofs of Theorem~\ref{thm:PM1} and Theorem~\ref{thm:PMR} rely on
the following two results for free boundary problems involving
degenerate and nonlipschitz parabolic equations.  We state them
independently because they could be interesting in themselves.

\begin{thm}\label{thm:FBP1}
	Let $x_{1}\leq x_{3}< x_{4}\leq x_{2}$, and let $c_{0}>0$ and
	$T>0$ be real numbers. Let $g:(0,c_{0})\to(0,+\infty)$ be a
	continuous function.
	Let $v:(x_{1},x_{2})\times[0,T)\to\re$ be a
	function such that
	\begin{enumerate}
		\renewcommand{\labelenumi}{(v\arabic{enumi})}
		\item  $v$ is continuous in $(x_{1},x_{2})\times[0,T)$;
	
		\item $v(x,t)\geq 0$ for every
		$(x,t)\in(x_{1},x_{2})\times[0,T)$;
	
		\item $v(x,0)> 0$ for every $x\in(x_{3},x_{4})$;

		\item  the partial derivative $v_{x}(x,t)$ exists for every 
		$(x,t)\in(x_{1},x_{2})\times(0,T)$;
	
		\item  setting
		$$\mathcal{P}:=\left\{(x,t)\in(x_{1},x_{2})\times(0,T):
		0<v(x,t)<c_{0}\right\},$$
		we have that $v\in C^{2,1}(\mathcal{P})$, and
		$$v_{t}(x,t)=g(v(x,t))v_{xx}(x,t)
		\quad\quad
		\forall (x,t)\in\mathcal{P}.$$
	\end{enumerate}
	
	Then $v(x,t)>0$ for every $(x,t)\in(x_{3},x_{4})\times[0,T)$.
\end{thm}

\begin{thm}\label{thm:FBP2}
	Let $r_{1}\leq r_{3}< r_{4}\leq r_{2}$, and let $c_{0}$, $c_{1}$,
	$T$, $G$, $A$ be positive real numbers.  Let
	$g:(0,c_{0})\to(0,+\infty)$ be a continuous function such that
	\begin{equation}
		\lim_{\sigma\to 0^{+}}\frac{g(\sigma)}{\sqrt{\sigma}}=G.
		\label{hp:g-lim}
	\end{equation}
	
	Let $f:(r_{1},r_{2})\times(0,T)\times[-c_{1},c_{1}]^{2}\to\re$ be 
	a function such that $f(r,t,0,0)=0$ uniformly in $(r,t)$, namely
	\begin{equation}
		\lim_{\sigma\to 0^{+}}\sup\left\{|f(r,t,p,q)|:
			(r,t,p,q)\in
			(r_{1},r_{2})\times(0,T)\times[-\sigma,\sigma]^{2}\right\}=0.
		\label{hp:FBP2-f}
	\end{equation}
	
	Let $v:(r_{1},r_{2})\times[0,T)\to\re$ be a
	function such that
	\begin{enumerate}
		\renewcommand{\labelenumi}{(v\arabic{enumi})}
		\item  $v$ is continuous in $(r_{1},r_{2})\times[0,T)$;
	
		\item $v(r,t)\geq 0$ for every
		$(r,t)\in(r_{1},r_{2})\times[0,T)$;
	
		\item $v(r,0)> 0$ for every $r\in(r_{3},r_{4})$;

		\item  the partial derivative $v_{r}(r,t)$ exists for every 
		$(r,t)\in(r_{1},r_{2})\times(0,T)$;
	
		\item  setting
		$$\mathcal{P}:=\left\{(r,t)\in(r_{1},r_{2})\times(0,T):
		0<v(r,t)<c_{0}\right\},$$
		we have that $v\in C^{2,1}(\mathcal{P})$, and
		\begin{equation}
			v_{t}\geq g(v)\left\{v_{rr}+f(r,t,v,v_{r})+A\right\}
			\quad\quad
			\forall (r,t)\in\mathcal{P}.
			\label{eq:FBP2-v}
		\end{equation}
	\end{enumerate}
	
	Finally, let $\mathcal{D}$ be the set defined as in (\ref{defn:D})
	with $k_{0}:=G\sqrt{A}$.
	
	Then $v(r,t)>0$ for every $(r,t)\in\mathcal{D}$.
\end{thm}

\setcounter{equation}{0}
\section{Proofs}\label{sec:proofs}

\subsection{Proof of Theorem~\ref{thm:PM1}}

We limit ourselves to prove that $u_{x}(x,t)<1$ for every
$(x,t)\in(x_{3},x_{4})\times[0,T)$. The proof of the symmetric
inequality $u_{x}(x,t)>-1$ is completely analogous.

\begin{figure}[htbp]
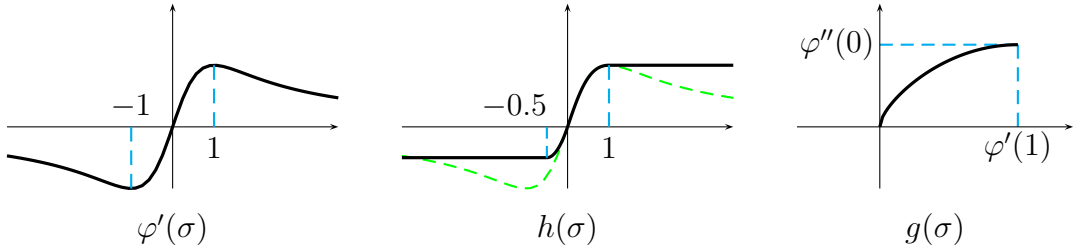

	\psset{unit=3ex}
	\centering
	\hfill
	\pspicture(-4,-3)(4,3)
	\psplot[linewidth=1.5\pslinewidth]{-4}{4}{x x 2 exp 1 add div 3 mul}
	\psline[linewidth=0.5\pslinewidth]{->}(-4,0)(4,0)
	\psline[linewidth=0.5\pslinewidth]{->}(0,-1.5)(0,3)
	\psline[linestyle=dashed, linecolor=cyan](1,0)(1,1.5)
	\psline[linestyle=dashed, linecolor=cyan](-1,0)(-1,-1.5)
	\rput(1,-0.5){1}
	\rput(-1,0.5){$-1$}
	\rput(0,-2.5){$\varphi'(\sigma)$}
	\endpspicture
	\hfill
	\pspicture(-4,-3)(4,3)
	\psplot[linecolor=green,linestyle=dashed]{-4}{4}{x x 2 exp 1 add div 3 mul}
	\psplot[linewidth=1.5\pslinewidth]{1}{4}{3 2 div}
	\psplot[linewidth=1.5\pslinewidth]{0}{1}{x x 2 exp 1 add div 3 mul}
	\psplot[linewidth=1.5\pslinewidth]{-0.5}{0}{x 0.5 add 2 exp 3 mul 3  
	4 div sub}
	\psplot[linewidth=1.5\pslinewidth]{-4}{-0.5}{-3 4 div}
	\psline[linewidth=0.5\pslinewidth]{->}(-4,0)(4,0)
	\psline[linewidth=0.5\pslinewidth]{->}(0,-1.5)(0,3)
	\psline[linestyle=dashed, linecolor=cyan](1,0)(1,1.5)
 	\psline[linestyle=dashed, linecolor=cyan](-0.5,0)(-0.5,-0.75)
	\rput(1,-0.5){1}
	\rput[r](-0.5,0.5){$-0.5$}
	\rput(0,-2.5){$h(\sigma)$}
	\endpspicture
	\hfill
	\psset{xunit=20ex}
	\pspicture(-0.3,-3)(0.7,3)
 	\psline[linestyle=dashed, linecolor=cyan](0.5,0)(0.5,2)
 	\psline[linestyle=dashed, linecolor=cyan](0,2)(0.5,2)
	\psplot[linewidth=1.5\pslinewidth]{0}{0.5}{x x x mul sub sqrt 2 x 
	x x mul sub mul add 2 mul}
	\psline[linewidth=0.5\pslinewidth]{->}(-0.3,0)(0.7,0)
	\psline[linewidth=0.5\pslinewidth]{->}(0,-1.5)(0,3)
 	\rput(0.5,-0.5){$\varphi'(1)$}
 	\rput[r](-0.02,2){$\varphi''(0)$}
	\rput(0.2,-2.5){$g(\sigma)$}
	\endpspicture
	\hfill\mbox{}
	\caption{Typical graph of functions $\varphi'$, $h$, and $g$}
	\label{fig:h}
\end{figure}

Let us introduce some notation.  Let us consider any function $h\in
C^{1}(\re)$ which is nondecreasing and such that
$h(\sigma)=\varphi'(\sigma)$ for every $\sigma\in[0,1]$,
$h(\sigma)=\varphi'(1)$ for every $\sigma\geq 1$, and $h(\sigma)$ is
constant for $\sigma\leq -1/2$.  Figure~\ref{fig:h} shows the typical
graph of such a function $h$. Note that condition (\ref{hp:phi=}) is
essential for the $C^{1}$ regularity of $h$.

The function $h$, as well as the function $\varphi'$, is invertible as
a function from $(0,1)$ to $(0,\varphi'(1))$. We can therefore define 
$g:(0,\varphi'(1))\to\re$ by setting
$$g(\sigma):=\varphi''\left(h^{-1}(\varphi'(1)-\sigma)\right)
\quad\quad
\forall\sigma\in(0,\varphi'(1)).$$

It is not difficult to see that the function $g$ is well defined,
positive, and continuous (but not Lipschitz continuous).  In the case
of the Perona-Malik equation with some computations one finds that
$g(\sigma)=(\sigma-\sigma^{2})^{1/2}+2(\sigma-\sigma^{2})$.  Its graph
is shown in Figure~\ref{fig:h}.

Let us finally set
$$v(x,t):=\varphi'(1)-h(u_{x}(x,t))
\quad\quad
\forall(x,t)\in(x_{1},x_{2})\times[0,T).$$

We claim that $g$ and $v$ satisfy the assumptions of
Theorem~\ref{thm:FBP1}. If we prove this claim, then we can conclude
that $v(x,t)>0$ in $(x_{3},x_{4})\times[0,T)$. This is equivalent to
say that $h(u_{x}(x,t))<\varphi'(1)$, which in turn is equivalent to
say that $u_{x}(x,t)<1$ in the same region.

So we only need to show that $v$ fulfils assumptions (v1) through (v5)
of Theorem~\ref{thm:FBP1}.

\paragraph{\textmd{\emph{Properties (v1) through (v3)}}}

The continuity of $v$ easily follows from the continuity of $h$ and
$u_{x}$.  Moreover $v(x,t)\geq 0$ in $(x_{1},x_{2})\times[0,T)$
because $h(\sigma)\leq\varphi'(1)$ for every $\sigma\in\re$.  Due to
(\ref{hp:PM1}) and the fact that $\varphi'(\sigma)<\varphi'(1)$ when
$\sigma<1$, we have that $v$ satisfies (v3).

\paragraph{\textmd{\emph{Property (v4)}}}

It is well known that $u$ is of class $C^{\infty}$ where
$|u_{x}(x,t)|\neq 1$ (because of the standard interior regularity
theory for parabolic equations).  Therefore the existence of the
partial derivative $v_{x}(x,t)$ is trivial when $|u_{x}(x,t)|\neq 1$.
The existence of $v_{x}(x,t)$ is trivial also when $u_{x}(x,t)=-1$
because $h(\sigma)$ is constant for $\sigma\leq 1/2$.

Let us consider now a point $(x_{0},t_{0})$ with
$u_{x}(x_{0},t_{0})=1$. We claim that in this point
$v_{x}(x_{0},t_{0})$ exists and is equal to 0. Let us assume that this
is not the case. Then there exists a sequence $\delta_{k}\to 0$ such
that
\begin{equation}
	\left|\frac{v(x_{0}+\delta_{k},t_{0})-
	v(x_{0},t_{0}))}{\delta_{k}}\right|\geq\nu>0
	\quad\quad
	\forall k\in\n.
	\label{eq:lim}
\end{equation}

Up to subsequences, we can always assume that either
$u_{x}(x_{0}+\delta_{k},t_{0})>1$ for every $k\in\n$, or
$u_{x}(x_{0}+\delta_{k},t_{0})\leq 1$ for every $k\in\n$.  In the
first case the fraction in (\ref{eq:lim}) is always 0, which is
incompatible with the condition stated therein.  In the second case
the fraction in (\ref{eq:lim}) can be rewritten as
$$-\frac{h(u_{x}(x_{0}+\delta_{k},t_{0}))-
h(u_{x}(x_{0},t_{0}))}{\delta_{k}}=-
\frac{\varphi'(u_{x}(x_{0}+\delta_{k},t_{0}))-
\varphi'(u_{x}(x_{0},t_{0}))}{\delta_{k}}.$$

When $\delta_{k}\to 0$, this quotient tends to
$(\varphi'(u_{x}))_{x}(x_{0},t_{0})$, and we know that this derivative
exists because $u$ is a solution of (\ref{eq:PM1}) of class $C^{1}$.
In order to find a contradiction, it is enough to show that this
derivative is equal to 0, and this is true because it is the
derivative in $x=x_{0}$ of the function
$x\to\varphi'(u_{x}(x,t_{0}))$, which attains its
maximum for $x=x_{0}$.

\paragraph{\textmd{\emph{Property (v5)}}}

Let us set
\begin{equation}
	\mathcal{P}:=\left\{(x,t)\in(x_{1},x_{2})\times(0,T):
	0<v(x,t)<\varphi'(1)\right\}.
	\label{defn:P1}
\end{equation}

From the properties of $h$ it follows that $0<u_{x}(x,t)<1$ in
$\mathcal{P}$, hence $u$ is of class $C^{\infty}$ in $\mathcal{P}$.
Since $h$ and $\varphi'$ coincide in $(0,1)$, we have therefore that
$$v_{t}=-h'(u_{x})u_{xt}=-h'(u_{x})\left(\varphi'(u_{x})\right)_{xx}=
-\varphi''(u_{x})(h(u_{x}))_{xx}=\varphi''(u_{x})(-h(u_{x}))_{xx}=$$
$$=\varphi''(u_{x})(\varphi'(1)-h(u_{x}))_{xx}=
\varphi''(u_{x})v_{xx}$$
in $\mathcal{P}$. Moreover in $\mathcal{P}$ we can express $u_{x}$ in 
function of $v$ as $u_{x}=h^{-1}\left(\varphi'(1)-v\right)$. In
conclusion
$$v_{t}=\varphi''(u_{x})v_{xx}=\varphi''\left(
h^{-1}\left(\varphi'(1)-v\right)\right)v_{xx}=g(v)v_{xx},$$
which proves that $v$ satisfies (v5).

\subsection{Proof of Theorem~\ref{thm:PMR}}

The argument is similar to the proof of Theorem~\ref{thm:PM1}. We
define $h$, $g$, and $v$ as in that proof, and we claim that the
assumptions of Theorem~\ref{thm:FBP2} are satisfied. 

\paragraph{\textmd{\emph{Properties of $g$}}}

As in the proof of Theorem~\ref{thm:PM1} it is quite easy to show that
$g:(0,\varphi'(1))\to(0,+\infty)$ is a continuous function.  It
remains to compute the limit in (\ref{hp:g-lim}).  Since we deal with
positive functions, we can square the numerator and the denominator.
Applying the variable change $\tau:=h^{-1}(\varphi'(1)-\sigma)$ and De
L'H\^{o}pital's rule, we have therefore that $$\lim_{\sigma\to
0^{+}}\frac{[g(\sigma)]^{2}}{\sigma}= \lim_{\sigma\to
0^{+}}\frac{[\varphi''\left(
h^{-1}\left(\varphi'(1)-\sigma\right)\right)]^{2}}{\sigma}=
\lim_{\tau\to 1^{-}}\frac{[\varphi''(\tau)]^{2}}{\varphi'(1)-
\varphi'(\tau)}=$$
$$=\lim_{\tau\to 1^{-}}\frac{2\varphi''(\tau)\varphi'''(\tau)}{-
\varphi''(\tau)}=-2\varphi'''(1)=2|\varphi'''(1)|,$$
which proves (\ref{hp:g-lim}) with $G:=\sqrt{2|\varphi'''(1)|}$.

\paragraph{\textmd{\emph{Properties of $v$}}}

The proof of (v1) through (v4) is analogous to the proof of the
corresponding properties in Theorem~\ref{thm:PM1}.  In order to prove
(v5), let us consider the set $\mathcal{P}$ defined in analogy with
(\ref{defn:P1}).  As in the previous case we have that
$0<u_{r}(r,t)<1$ in this set, hence $v$ is regular and satisfies
\begin{equation}
	v_{t}=-h'(u_{r})u_{rt}=-\varphi''(u_{r})\left[
	\left(\varphi'(u_{r})\right)_{r}+
	\frac{\varphi'(u_{r})}{r}\right]_{r}.
	\label{eq:PMR-vt}
\end{equation}

Once again $\varphi''(u_{r})=\varphi''\left(
h^{-1}\left(\varphi'(1)-v\right)\right)=g(v)$.  Moreover 
$$\left[ \left(\varphi'(u_{r})\right)_{r}+
\frac{\varphi'(u_{r})}{r}\right]_{r}=
\left(\varphi'(u_{r})\right)_{rr}+
\frac{\left(\varphi'(u_{r})\right)_{r}}{r}-
\frac{\varphi'(u_{r})}{r^{2}}=
-v_{rr}-\frac{v_{r}}{r}-\frac{\varphi'(1)}{r^{2}}+\frac{v}{r^{2}}.$$

Plugging these identities into (\ref{eq:PMR-vt}) we obtain that
$$v_{t}=g(v)\left\{v_{rr}+\frac{v_{r}}{r}-
\frac{v}{r^{2}}+\frac{\varphi'(1)}{r^{2}}\right\}\geq
g(v)\left\{v_{rr}+f(r,t,v,v_{r})+A\right\},$$
where
$$f(r,t,p,q):=\frac{q}{r}-\frac{p}{r^{2}},
\hspace{3em}
A:=\frac{\varphi'(1)}{r_{2}^{2}}.$$

We have thus proved that $v$ satisfies the differential inequality
(\ref{eq:FBP2-v}) in $\mathcal{P}$ with a constant $A>0$, and a
function $f(r,t,p,q)$ satisfying (\ref{hp:FBP2-f}).

\paragraph{\textmd{\emph{Conclusion}}}

From Theorem~\ref{thm:FBP2} we deduce that $v$ is positive in the
region $\mathcal{D}$ defined according to (\ref{defn:D}) with
$k_{0}=G\sqrt{A}=r_{2}^{-1}\sqrt{2\varphi'(1)|\varphi'''(1)|}$, which
in turn implies that $u_{r}<1$ in the same region.  The proof of the
symmetric inequality $u_{r}>-1$ is completely analogous.

\subsection{Proof of Theorem~\ref{thm:no-ex}}

Let us assume that a solution exists with 
$$T>T_{0}:=\frac{r_{2}(r_{2}-r_{1})}{\sqrt{2\varphi'(1)|\varphi'''(1)|}}.$$

Let us set
$$M(t):=\max\left\{u_{r}(r,t):r\in[r_{3},r_{5}]\right\}
\quad\quad
\forall t\in[0,T).$$

Due to the first two inequalities in (\ref{hp:i-data}) we know that
the subcritical region is nonempty at time $t=0$.  Applying
Theorem~\ref{thm:PMR} we have therefore that the subcritical region
expands, and coincides with the whole interval $(r_{1},r_{2})$ as soon
as
$$t>\frac{r_{2}-r_{1}}{k_{0}}=T_{0}.$$

In particular this means that
\begin{equation}
	M(T_{0})\leq 1.
	\label{eq:M(t0)}
\end{equation}

On the other hand we claim that
\begin{equation}
	M(t)\geq M(0)-\frac{\varphi'(1)}{r_{1}^{2}}t
	\quad\quad
	\forall t\in[0,T_{0}].
	\label{est:M(t)}
\end{equation}

If we prove this claim, then setting $t=T_{0}$ and exploiting
the last inequality in (\ref{hp:i-data}), we find that
\begin{eqnarray*}
	M(T_{0}) & \geq & M(0)-\frac{\varphi'(1)}{r_{1}^{2}}T_{0}  \\
	& > & 1+\frac{r_{2}(r_{2}-r_{1})}{r_{1}^{2}}
	 \sqrt{\frac{\varphi'(1)}{2|\varphi'''(1)|}}-
	 \frac{\varphi'(1)}{r_{1}^{2}}\cdot
	 \frac{r_{2}(r_{2}-r_{1})}{\sqrt{2\varphi'(1)|\varphi'''(1)|}}\\
	 & = & 1,
\end{eqnarray*}
which contradicts (\ref{eq:M(t0)}).

\paragraph{\textmd{\emph{Proof of (\ref{est:M(t)})}}}

The argument is similar to the usual comparison principles.  Setting
for simplicity $v(r,t):=u_{r}(r,t)$, we have that $v$ is a solution of
\begin{equation}
	v_{t}=\varphi''(v)v_{rr}+\varphi'''(v)v_{r}^{2}+
	\frac{\varphi''(v)}{r}v_{r}-\frac{\varphi'(1)}{r^{2}}
	\label{eq:v}
\end{equation}
in the subset of $(r_{1},r_{2})\times[0,T)$ where $|v|\neq 1$. Let us set
\begin{equation}
	w(t):=M(0)-\ep-\frac{\varphi'(1)}{r_{1}^{2}}t
	\quad\quad
	\forall t\in[0,T_{0}],
	\label{defn:w}
\end{equation}
where $\ep>0$ is small enough so that $w(T_{0})>1$, hence $w(t)>1$ for
every $t\in[0,T_{0}]$. We claim that 
\begin{equation}
	M(t)\geq w(t)
	\quad\quad
	\forall t\in[0,T_{0}],
	\label{eq:M>w}
\end{equation}
from which (\ref{est:M(t)}) follows by letting $\ep\to 0^{+}$. Let us 
prove (\ref{eq:M>w}) by contradiction. Let us assume that $M(t)<w(t)$ 
for some $t\in[0,T_{0}]$, and let us set
$$t_{0}:=\inf\left\{t\in[0,T_{0}]:v(r,t)<w(t)\quad\forall
r\in[r_{3},r_{5}]\right\}.$$

Since $M(0)>w(0)$, we have that $t_{0}>0$.  Moreover, due to the
continuity of $v$ and $w$, there exists $r_{0}\in[r_{3},r_{5}]$ such
that $v(r_{0},t_{0})=w(t_{0})$, and 
$$v(r,t_{0})-w(t_{0})\leq 0
\quad\quad
\forall r\in[r_{3},r_{5}].$$

Since subcritical regions don't shrink, we have that
$v(r_{3},t_{0})<1$ and $v(r_{5},t_{0})<1$, while $w(t_{0})>1$.  This
shows in particular that $r_{0}\neq r_{3}$ and $r_{0}\neq r_{5}$.  Now
we know that $r_{0}$ is a maximum point for the function $r\to
v(r,t_{0})-w(t_{0})$, and $r_{0}$ is contained in the open interval
$(r_{3},r_{5})$, hence
\begin{equation}
	v_{r}(r_{0},t_{0})=0\quad\mbox{ and }\quad
	v_{rr}(r_{0},t_{0})\leq 0.
	\label{no-ex:vr-vrr}
\end{equation}

Let us consider now time derivatives.  Since
$v(r_{0},t_{0})=w(t_{0})>1$, we can use (\ref{eq:v}).  Exploiting also
(\ref{defn:w}) and (\ref{no-ex:vr-vrr}) we obtain that
\begin{eqnarray*}
	v_{t}(r_{0},t_{0})-w_{t}(t_{0}) & = &
	\varphi''(v(r_{0},t_{0}))v_{rr}(r_{0},t_{0})-
	\frac{\varphi'(1)}{r_{0}^{2}}+\frac{\varphi'(1)}{r_{1}^{2}}\\
	 & > & \varphi''(v(r_{0},t_{0}))v_{rr}(r_{0},t_{0}).
\end{eqnarray*}

Since $\varphi''(v(r_{0},t_{0}))\leq 0$, we can conclude that
$$v_{t}(r_{0},t_{0})-w_{t}(t_{0})
> \varphi''(v(r_{0},t_{0}))v_{rr}(r_{0},t_{0})\geq 0.$$

This implies that $v(r_{0},t)-w(t)>0$ for every $t$ in a suitable 
right-hand neighborhood of $t_{0}$, which contradicts the definition
of $t_{0}$.

\subsection{Proof of Theorem~\ref{thm:example}}

Let us set
\begin{equation}
	k_{1}:=n,
	\quad\quad
	k_{2}:=1,
	\quad\quad
	h_{1}:=n^{3},
	\quad\quad
	h_{2}:=-n^{2}.
	\label{defn:hk}
\end{equation}

We claim that the conclusions of statement~(1) and statement~(2) hold
true provided that $n$ is large enough.

\paragraph{\textmd{\emph{Statement (1)}}}

Let us set for simplicity $v_{0}(x,y):=|\nabla u_{0}(x,y)|^{2}$. Let
us assume that
\begin{equation}
	v_{0y}(0,0)>0.
	\label{eq:dini}
\end{equation}
 
Then the implicit function theorem implies that the set $v_{0}(x,y)>1$
can be represented, in a neighborhood of $(0,0)$, as $y>g(x)$, where
$g$ is a suitable function defined in a neighborhood of $x=0$. Such a 
function satisfies 
$$g(0)=0,
\hspace{3em}
g'(0)=-\frac{v_{0x}(0,0)}{v_{0y}(0,0)},$$
\begin{equation}
	g''(0)=-\frac{1}{v_{0y}^{3}}\left\{
	v_{0x}^{2}v_{0yy}+v_{0y}^{2}v_{0xx}
	-2v_{0x}v_{0y}v_{0xy}\right\},
	\label{eq:g''(0)}
\end{equation}
where in (\ref{eq:g''(0)}) all partial derivatives of $v_{0}$ are
computed in $(0,0)$.  In particular $g$ is convex in a neighborhood of
$0$ if the right-hand side of (\ref{eq:g''(0)}) is positive.  

From (\ref{hp:taylor}) we have that, up to higher order terms,
$$u_{0x}(x,y)=\frac{\sqrt{2}}{2}+2k_{1}x+3h_{1}x^{2},
\quad\quad
u_{0y}(x,y)=\frac{\sqrt{2}}{2}+2k_{2}y+3h_{2}y^{2},$$
hence
$$v_{0}(x,y)=1+2\sqrt{2}(k_{1}x+k_{2}y)+
\left(4k_{1}^{2}+3\sqrt{2}h_{1}\right)x^{2}+
\left(4k_{2}^{2}+3\sqrt{2}h_{2}\right)y^{2}.$$

All the derivatives appearing in (\ref{eq:dini}) and (\ref{eq:g''(0)})
can be easily computed.  It follows that condition (\ref{eq:dini}) is
equivalent to $k_{2}>0$, while $g''(0)>0$ if and only if
$$8k_{1}^{2}k_{2}^{2}+3\sqrt{2}
\left(k_{1}^{2}h_{2}+k_{2}^{2}h_{1}\right)<0.$$

Both conditions are satisfied if the values of the parameters are
given by (\ref{defn:hk}) and $n$ is large enough.

\paragraph{\textmd{\emph{Statement (2)}}}

Let us set for simplicity $v(x,y,t):=|\nabla u(x,y,t)|^{2}$. Thesis is
proved if we show that
$$v_{t}=2u_{x}u_{tx}+2u_{y}u_{ty}>0$$
in the point $(x,y,t)=(0,0,0)$. We can therefore deduce the value of
$v_{t}(0,0,0)$ from the Taylor expansion of $u_{0}(x,y)$.

In order to compute $u_{t}$, we recall that
$u_{t}=\Psi_{1x}+\Psi_{2y}$, where
$$\Psi_{1}:=\varphi'\left(\left(u_{x}^{2}+u_{y}^{2}\right)^{1/2}\right)
\frac{u_{x}}{\left(u_{x}^{2}+u_{y}^{2}\right)^{1/2}},
\quad\quad
\Psi_{2}:=\varphi'\left(\left(u_{x}^{2}+u_{y}^{2}\right)^{1/2}\right)
\frac{u_{y}}{\left(u_{x}^{2}+u_{y}^{2}\right)^{1/2}}.$$

With some computations we obtain that, up to higher order terms,
$$\left(u_{x}^{2}+u_{y}^{2}\right)^{1/2}=
1+\sqrt{2}(k_{1}x+k_{2}y)+\left[
k_{1}^{2}+\frac{3\sqrt{2}}{2}h_{1}\right]x^{2} +\left[
k_{2}^{2}+\frac{3\sqrt{2}}{2}h_{2}\right]y^{2}-2k_{1}k_{2}xy,$$
hence
$$\varphi'\left(\left(u_{x}^{2}+u_{y}^{2}\right)^{1/2}\right)=
\varphi'(1)+ (k_{1}x+k_{2}y)^{2}\varphi'''(1),$$
and therefore
\begin{eqnarray*}
	\Psi_{1} & = &\frac{\sqrt{2}}{2}\varphi'(1)+\varphi'(1)\left(
	k_{1}x-k_{2}y\right)+\frac{1}{2}\left[
	3\varphi'(1)\left(h_{1}-\sqrt{2}\,k_{1}^{2}\right)+
	\sqrt{2}\,\varphi'''(1)k_{1}^{2}\right]x^{2}  \\
	\noalign{\vspace{1ex}}
	 &  & +\frac{1}{2}\left[\varphi'(1)\left(\sqrt{2}\,k_{2}^{2}-3h_{2}\right)
	+\sqrt{2}\varphi'''(1)k_{2}^{2}\right]y^{2}+
	\sqrt{2}\left(\varphi'(1)+\varphi'''(1)\strut\right)k_{1}k_{2}xy.
\end{eqnarray*}

The expression for $\Psi_{2}$ is just the symmetric one.  It follows
that, up to higher order terms,
\begin{eqnarray*}
	u_{t}(x,y,0) & = & (k_{1}+k_{2})\varphi'(1) \\
	 & & +
	\left\{\left(3h_{1}+\sqrt{2}\,k_{1}(k_{2}-3k_{1})\right)\varphi'(1)+
	\sqrt{2}\varphi'''(1)k_{1}(k_{1}+k_{2})\right\}x \\
	 &  & +\left\{\left(3h_{2}+\sqrt{2}\,k_{2}(k_{1}-3k_{2})\right)\varphi'(1)+
	\sqrt{2}\varphi'''(1)k_{2}(k_{1}+k_{2})\right\}y,
\end{eqnarray*}
hence
$$v_{t}(0,0,0)=\varphi'(1)\left\{3\sqrt{2}(h_{1}+h_{2})+
4k_{1}k_{2}-6k_{1}^{2}-6k_{2}^{2}\right\}+
2\varphi'''(1)(k_{1}+k_{2})^{2}.$$

From this expression it is easy to see that $v_{t}(0,0,0)>0$ if the
values of the parameters are given by (\ref{defn:hk}) and $n$ is large
enough.

\subsection{Proof of Theorem~\ref{thm:FBP1}}

Let $x_{\star}\in(x_{3},x_{4})$ be any point. We have to prove that
\begin{equation}
	v(x_{\star},t)>0
	\quad\quad
	\forall t\in[0,T).
	\label{eq:th-FBP1}
\end{equation}

To this end we fix some notation. First of all we choose real numbers 
$x_{5}$ and $x_{6}$ such that $x_{3}<x_{5}<x_{\star}<x_{6}<x_{4}$.
Then we consider the functions 
$$\psi(x):=(x-x_{5})(x_{6}-x),
\quad\quad\quad
w(x,t):=e^{-\lambda t}\left(\delta^{2}\psi(x)+
\delta\psi^{2}(x)\right),$$
where $\lambda$ and $\delta$ are positive parameters. We claim that
when $\lambda$ is large enough and $\delta$ is small enough we have
that
\begin{equation}
	v(x,t)\geq w(x,t)
	\quad\quad
	\forall (x,t)\in[x_{5},x_{6}]\times[0,T).
	\label{FBP1:v>w}
\end{equation}

Since $w$ is positive in $(x_{5},x_{6})\times[0,T)$, and
$x_{\star}\in(x_{5},x_{6})$, this is enough to prove
(\ref{eq:th-FBP1}). In order to prove (\ref{FBP1:v>w}) we first
establish some properties of $w$.

\paragraph{\textmd{\emph{Properties of $w$}}}

Let us show that $w$ fulfils the following properties.
\begin{enumerate}
	\renewcommand{\labelenumi}{(w\arabic{enumi})}
	\item  $w\in C^{\infty}\left([x_{5},x_{6}]\times[0,+\infty)\right)$;

	\item $w(x,t)>0$ for every
	$(x,t)\in(x_{5},x_{6})\times[0,+\infty)$;

	\item $w_{x}(x_{5},t)=\delta^{2}e^{-\lambda t}\psi'(x_{5})>0$ for
	every $t\geq 0$;

	\item $w_{x}(x_{6},t)=\delta^{2}e^{-\lambda t}\psi'(x_{6})<0$ for
	every $t\geq 0$;

	\item  if $\delta$ is small enough we have that $w(x,0)<v(x,0)$
	for every $x\in[x_{5},x_{6}]$;

	\item if $\delta$ is small enough we have that $w(x,t)<c_{0}$ for
	every $(x,t)\in(x_{5},x_{6})\times[0,+\infty)$;

	\item  if $\delta$ is small enough and $\lambda$ is large enough, 
	then $w$ satisfies
	\begin{equation}
		w_{t}(x,t)<g(w(x,t))w_{xx}(x,t)
		\quad\quad
		\forall(x,t)\in(x_{5},x_{6})\times[0,+\infty).
		\label{FBP1:eq-w}
	\end{equation}
\end{enumerate}

Properties (w1) through (w4) easily follow from the definition of $w$ 
and $\psi$. Property (w5) follows from the fact that the infimum of
$v(x,0)$ for $x\in[x_{5},x_{6}]$ is strictly positive due to (v3).
Property (w6) is almost trivial. In order to prove (w7) we recall that
$\psi''(x)=-2$, hence
$$w_{t}(x,t)=-\lambda e^{-\lambda t}\left(\delta^{2}\psi(x)+
\delta\psi^{2}(x)\right),
\quad
w_{xx}(x,t)=e^{-\lambda
t}\left(2\delta[\psi'(x)]^{2}-4\delta\psi(x)-2\delta^{2}\right).$$

Plugging these identities in (\ref{FBP1:eq-w}), we are left to prove
that
\begin{equation}
	-\lambda\left(\delta\psi+\psi^{2}\right)<g(w)\left(
	2[\psi']^{2}-4\psi-2\delta\right)
	\quad\quad
	\forall(x,t)\in(x_{5},x_{6})\times[0,+\infty).
	\label{FBP1:eq-w-true}
\end{equation}

To this end we fix once for all two real numbers $x_{7}$ and $x_{8}$
such that $x_{5}<x_{7}<x_{8}<x_{6}$, and
$$\inf\left\{2[\psi'(x)]^{2}-4\psi(x):
x\in[x_{5},x_{7}]\cup[x_{8},x_{6}]\right\}>0.$$

This is possible because in the endpoints of the interval
$[x_{5},x_{6}]$ one has that $\psi=0$ and $\psi'\neq 0$.  Now we
distinguish two cases.  When $x\in[x_{5},x_{7}]\cup[x_{8},x_{6}]$ the
left-hand side of (\ref{FBP1:eq-w-true}) is negative, while the
left-hand side is positive provided that $\delta$ is small enough,
independently on $x$.  When $x\in[x_{7},x_{8}]$ the right-hand side
may be negative, but also the left-hand side is strictly negative
because in this interval $\psi$ is bounded from below by a positive
constant.  In other words, in $[x_{7},x_{8}]\times[0,+\infty)$
inequality (\ref{FBP1:eq-w-true}) holds true if we choose
$$\lambda>\sup\left\{\frac{g(w(x,t))\left(4\psi(x)+
2\delta\right)}{\delta\psi(x)+\psi^{2}(x)}
:(x,t)\in[x_{7},x_{8}]\times[0,+\infty)\right\}.$$

We point out that the supremum is finite.  This completes
the proof of (\ref{FBP1:eq-w-true}).

\paragraph{\textmd{\emph{Proof of (\ref{FBP1:v>w})}}}

Let us choose positive values of $\delta$ and $\lambda$ in such a way
that $w$ satisfies (w1) through (w7).  Now we argue more or less as in
the proof of the classical comparison results.  Let us assume that
(\ref{FBP1:v>w}) is not true, and let us set
$$t_{0}:=\inf\left\{t\in[0,T):\exists x\in[x_{5},x_{6}]
\mbox{ such that }v(x,t)-w(x,t)<0\right\}.$$

From (w5) we have that $t_{0}>0$. Moreover, from the definition of
$t_{0}$ it follows that
\begin{equation}
	v(x,t)-w(x,t)\geq 0
	\quad\quad
	\forall(x,t)\in[x_{5},x_{6}]\times [0,t_{0}].
	\label{FBP1:v-w-0}
\end{equation}

Finally, from the continuity of $v$ and $w$ we deduce that there
exists $x_{0}\in[x_{5},x_{6}]$ such that
$v(x_{0},t_{0})-w(x_{0},t_{0})=0$.

We claim that $x_{0}\neq x_{5}$. Indeed let us assume by
contradiction that $x_{0}=x_{5}$. Then $w(x_{0},t_{0})=0$, hence also 
$v(x_{0},t_{0})=0$. By (v2) it follows that $x_{0}$ is a minimum 
point for the function $x\to v(x,t_{0})$. By (v4) we have therefore
that $v_{x}(x_{0},t_{0})=0$. Keeping (w3) into account, we deduce that
\begin{equation}
	(v-w)_{x}(x_{0},t_{0})=v_{x}(x_{0},t_{0})-w_{x}(x_{0},t_{0})<0.
	\label{FBP1:v-w-x}
\end{equation}

On the other hand, from (\ref{FBP1:v-w-0}) we know also that
$x_{0}=x_{5}$ is a minimum point for the function $x\to
v(x,t_{0})-w(x,t_{0})$ restricted to the interval $[x_{5},x_{6}]$.
Since the minimum point is the left-hand endpoint of the interval, we
deduce that $(v-w)_{x}(x_{0},t_{0})\geq 0$, which contradicts
(\ref{FBP1:v-w-x}).

In a symmetric way we prove that $x_{0}\neq x_{6}$.  So we are left
with the case in which $t_{0}>0$ and $x_{0}\in(x_{5},x_{6})$.  In this
case by (w2) and (w6) we have that
$0<v(x_{0},t_{0})=w(x_{0},t_{0})<c_{0}$, hence both $v$ and $w$ are
smooth in a neighborhood of this point and fulfil (v5) and (w7),
respectively.  In particular, since $x_{0}$ is always a minimum point
of the function $x\to v(x,t_{0})-w(x,t_{0})$, and now $x_{0}$ is in
the interior of the interval $(x_{5},x_{6})$, we have that
\begin{equation}
	v_{x}(x_{0},t_{0})= w_{x}(x_{0},t_{0})\quad\mbox{ and }\quad
	v_{xx}(x_{0},t_{0})\geq w_{xx}(x_{0},t_{0}).
	\label{FBP1:vxx-wxx}
\end{equation}

Let us consider now the time derivatives.  On the one hand, from
(\ref{FBP1:v-w-0}) we deduce that $(v-w)_{t}(x_{0},t_{0})\leq 0$.  On the
other hand, from (v5), (w7), and (\ref{FBP1:vxx-wxx}) we have that
\begin{eqnarray*}
	(v-w)_{t}(x_{0},t_{0}) & > &
	g(v(x_{0},t_{0}))v_{xx}(x_{0},t_{0})-
	g(w(x_{0},t_{0}))w_{xx}(x_{0},t_{0}) \\
	 & = & g(w(x_{0},t_{0}))\left(
	 v_{xx}(x_{0},t_{0})-w_{xx}(x_{0},t_{0})\right) \\
	 & \geq & 0.
\end{eqnarray*}

This rules out the last case and completes the proof of
(\ref{FBP1:v>w}).

\subsection{Proof of Theorem~\ref{thm:FBP2}}

The strategy is similar to the proof of Theorem~\ref{thm:FBP1}.  The
main difference is that in this case we have to cope with moving
domains.

Let $(r_{\star},t_{\star})$ be any point of $\mathcal{D}$.  We have to
prove that
\begin{equation}
	v(r_{\star},t_{\star})>0.
	\label{eq:th-FBP2}
\end{equation}

To this end we fix some notation.  First of all it is not difficult to
see that there exists real numbers $r_{5}$, $r_{6}$, $k$ such that 
$$|k|<G\sqrt{A},
\quad\quad\quad
r_{1}<r_{5}+kt_{\star}<r_{\star}<r_{6}+kt_{\star}<r_{2}.$$

Then we consider the set
$$\mathcal{D}_{\star}:=\left\{(r,t)\in(r_{1},r_{2})\times[0,T): t\leq
t_{\star},\ r_{5}+kt\leq r\leq r_{6}+kt\right\}\subseteq\mathcal{D}.$$

\begin{figure}[htbp]
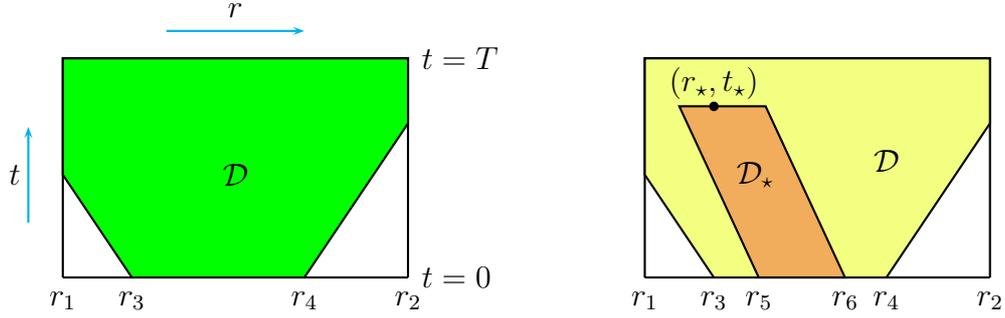

\psset{unit=5ex}
\centering
\hfill
\pspicture(-1,-0.5)(6,4.2)
\psline(0,0)(5,0)
\psline(0,0)(0,3.2)(5,3.2)(5,0)
\pspolygon[fillstyle=solid,fillcolor=green](1,0)(0,1.5)(0,3.2)(5,3.2)(5,2.25)(3.5,0)(1,0)
\psline[linecolor=cyan]{->}(-0.5,0.8)(-0.5,2.2)
\psline[linecolor=cyan]{->}(1.5,3.6)(3.5,3.6)
\rput[t](0,-0.2){$r_{1}$}
\rput[t](1,-0.2){$r_{3}$}
\rput[t](3.5,-0.2){$r_{4}$}
\rput[t](5,-0.2){$r_{2}$}
\rput[l](5.2,0){$t=0$}
\rput[l](5.2,3.2){$t=T$}
\rput(-0.7,1.5){$t$}
\rput(2.5,3.9){$r$}
\rput(2.5,1.5){$\mathcal{D}$}
\endpspicture
\hfill
\pspicture(-1,-0.5)(5,4.2)
\newrgbcolor{fuoco}{0.95 0.68 0.36}
\newrgbcolor{giallino}{0.95 1 0.5}
\psline(0,0)(5,0)
\psline(0,0)(0,3.2)(5,3.2)(5,0)
\pspolygon[fillstyle=solid,fillcolor=giallino](1,0)(0,1.5)(0,3.2)(5,3.2)(5,2.25)(3.5,0)(1,0)
\pspolygon[fillstyle=solid,fillcolor=fuoco](1.65,0)(0.5,2.5)(1.75,2.5)(2.9,0)
\psdots(1,2.5)
\rput[t](0,-0.2){$r_{1}$}
\rput[t](1,-0.2){$r_{3}$}
\rput[t](3.5,-0.2){$r_{4}$}
\rput[t](1.65,-0.2){$r_{5}$}
\rput[t](2.9,-0.2){$r_{6}$}
\rput[t](5,-0.2){$r_{2}$}
\rput(1.6,1.5){$\mathcal{D}_{\star}$}
\rput(3.5,1.7){$\mathcal{D}$}
\rput[b](1,2.6){$(r_{\star},t_{\star})$}
\endpspicture
\hfill\mbox{}
	\caption{the sets $\mathcal{D}$ and $\mathcal{D}_{\star}$}
	\label{fig:D-star}
\end{figure}

We refer to Figure~\ref{fig:D-star} for a representation of the set
$\mathcal{D}_{\star}$ (corresponding in that case to some $k<0$) and
its relation with $(r_{\star},t_{\star})$ and $\mathcal{D}$.  Note
that the slope of the slanted lines bounding $\mathcal{D}_{\star}$ is
larger than the slope of the slanted lines limiting $\mathcal{D}$.
This is just because $|k|<G\sqrt{A}$.

Due to this inequality, there exists $\ep_{0}\in(0,\min\{1,G,A/2\})$
such that
\begin{equation}
	|k|<(1-\ep_{0})(G-\ep_{0})\sqrt{A-2\ep_{0}}.
	\label{defn:ep0}
\end{equation}

From now on, $\ep_{0}$ is a fixed positive constant.  Due to
(\ref{hp:g-lim}) there exists also $c_{2}\in(0,c_{0})$ such that
\begin{equation}
	g(\sigma)\geq(G-\ep_{0})\sqrt{\sigma} \quad\quad
	\forall\sigma\in(0,c_{2}).
	\label{est:gep}
\end{equation}

Finally we consider the functions 
$$\psi(r):=(r-r_{5})(r_{6}-r),
\quad\quad\quad
w(r,t):=\delta^{3}\psi(r-kt)+\delta[\psi(r-kt)]^{3/2},$$
where $\delta$ is a positive parameter.  We claim that when $\delta$
is small enough we have that
\begin{equation}
	v(r,t)\geq w(r,t)
	\quad\quad
	\forall (r,t)\in\mathcal{D}_{\star}.
	\label{FBP2:v>w}
\end{equation}

This inequality, applied with $(r,t)=(r_{\star},t_{\star})$, implies 
(\ref{eq:th-FBP2}). In order to prove (\ref{FBP2:v>w}) we first
establish some properties of $w$.

\paragraph{\textmd{\emph{Properties of $w$}}}

Let $\mbox{Int}(\mathcal{D}_{\star})$ denote the set of points
$(r,t)\in\mathcal{D}_{\star}$ with $r_{5}+kt<r<r_{6}+kt$
Let us show that $w$ fulfils the following properties.
\begin{enumerate}
	\renewcommand{\labelenumi}{(w\arabic{enumi})}
	\item  $w\in C^{1}\left(\mathcal{D}_{\star}\right)\cap
	C^{\infty}\left(\mbox{Int}(\mathcal{D}_{\star})\right)$;

	\item $w(r,t)>0$ for every
	$(r,t)\in\mbox{Int}(\mathcal{D}_{\star})$;

	\item $w_{r}(r_{5}+kt,t)=\delta^{3}\psi'(r_{5})>0$ for every
	$t\in[0,t_{\star}]$;

	\item $w_{r}(r_{6}+kt,t)=\delta^{3}\psi'(r_{6})<0$ for every
	$t\in[0,t_{\star}]$;

	\item  if $\delta$ is small enough we have that $w(r,0)<v(r,0)$
	for every $r\in[r_{5},r_{6}]$;

	\item if $\delta$ is small enough we have that $w(r,t)<c_{2}$ for
	every $(r,t)\in\mathcal{D}_{\star}$;

	\item if $\delta$ is small enough, then $w$ satisfies
	\begin{equation}
		w_{t}<g(w)\left\{w_{rr}+f(r,t,w,w_{r})+A\right\}
		\quad\quad
		\forall(r,t)\in\mbox{Int}(\mathcal{D}_{\star}).
		\label{FBP2:eq-w}
	\end{equation}
\end{enumerate}

Properties (w1) through (w4) easily follow from the definition of $w$ 
and $\psi$. Property (w5) follows from the fact that the infimum of
$v(r,0)$ for $r\in[r_{5},r_{6}]$ is strictly positive due to (v3).
Property (w6) is almost trivial. In order to prove (w7) we recall that
$\psi''(r)=-2$, hence (for simplicity we set $y:=r-kt$, and we observe
that $y\in[r_{5},r_{6}]$)
$$w_{t}(r,t)=-k\delta^{3}\psi'(y)-\frac{3}{2}k\delta
\left[\psi(y)\right]^{1/2}\psi'(y),$$
$$w_{r}(r,t)=\delta^{3}\psi'(y)+\frac{3}{2}\delta
\left[\psi(y)\right]^{1/2}\psi'(y),$$
$$w_{rr}(r,t)=-2\delta^{3}-3\delta
\left[\psi(y)\right]^{1/2}+\frac{3}{4}\delta
\left[\psi(y)\right]^{-1/2}\left[\psi'(y)\right]^{2}.$$

When $\delta\to 0^{+}$ we have that $w$ and $w_{r}$ tend to zero
uniformly in $\mathcal{D}_{\star}$. Thanks to (\ref{hp:FBP2-f}) we
have therefore that
$$\left|f(r,t,w(r,t),w_{r}(r,t))\right|\leq\ep_{0}
\quad\quad
\forall(r,t)\in\mathcal{D}_{\star},$$
provided that $\delta$ is small enough. In an analogous way we have
also that
$$w_{rr}(r,t)>\frac{3}{4}\delta
\left[\psi(y)\right]^{-1/2}\left[\psi'(y)\right]^{2}-\ep_{0}$$
provided that $\delta$ is small enough. It follows that
$$w_{rr}+f(r,t,w,w_{r})+A> A-2\ep_{0}+
\frac{3}{4}\delta
\left[\psi(y)\right]^{-1/2}\left[\psi'(y)\right]^{2}>0$$
in $\mbox{Int}(\mathcal{D}_{\star})$.  Moreover, from (\ref{est:gep}),
(w2) and (w6) we have that
$$g(w(r,t))\geq(G-\ep_{0})\sqrt{w(r,t)}
\quad\quad
\forall(r,t)\in\mbox{Int}(\mathcal{D}_{\star}),$$
and in conclusion
$$g(w)\left\{w_{rr}+f(r,t,w,w_{r})+A\right\}>
(G-\ep_{0})\sqrt{w}\left\{A-2\ep_{0}+
\frac{3}{4}\delta
\frac{\left[\psi'(y)\right]^{2}}{\left[\psi(y)\right]^{1/2}}
\right\}.$$

Therefore inequality (\ref{FBP2:eq-w}) is proved if we show that
\begin{equation}
	\delta^{3}|k\psi'|+\frac{3}{2}\delta|k\psi'|\sqrt{\psi}\leq
	(G-\ep_{0})\sqrt{\delta^{3}\psi+\delta[\psi]^{3/2}}
	\left\{A-2\ep_{0}+
	\frac{3}{4}\delta
	\frac{[\psi']^{2}}{\sqrt{\psi}}
	\right\},
	\label{FBP2:eq-w-true}
\end{equation}
where the argument of $\psi$ and $\psi'$ is any $y\in(r_{5},r_{6})$.
Let us consider the right-hand side of (\ref{FBP2:eq-w-true})
multiplied by $(1-\ep_{0})$.  Applying the inequality between the
arithmetic mean and the geometric mean, and recalling (\ref{defn:ep0}),
we obtain that
\begin{eqnarray*}
	(1-\ep_{0})\cdot(\mbox{right-hand side}) & \geq &
	(1-\ep_{0})(G-\ep_{0})\sqrt{\delta[\psi]^{3/2}}\cdot
	\left\{A-2\ep_{0}+ \frac{3}{4}\delta
	\frac{[\psi']^{2}}{\sqrt{\psi}} \right\}\\
	 & \geq & (1-\ep_{0})(G-\ep_{0})\sqrt{\delta[\psi]^{3/2}}\cdot
	2\left[(A-2\ep_{0})\cdot\frac{3}{4}\delta
	\frac{[\psi']^{2}}{\sqrt{\psi}}\right]^{1/2}  \\
	 & = & (1-\ep_{0})(G-\ep_{0})
	\sqrt{A-2\ep_{0}}\cdot\delta|\psi'|\sqrt{\psi}\cdot\sqrt{3}  \\
	 & \geq & \frac{3}{2}\delta|k\psi'|\sqrt{\psi}.
\end{eqnarray*}

In order to prove (\ref{FBP2:eq-w-true}) it is therefore enough to
show that
$$\delta^{3}|k\psi'|\leq\ep_{0}
(G-\ep_{0})\sqrt{\delta^{3}\psi+\delta[\psi]^{3/2}} \left\{A-2\ep_{0}+
\frac{3}{4}\delta \frac{[\psi']^{2}}{\sqrt{\psi}} \right\},$$
which in turn is true if we show that
$$\delta^{3/2}|k\psi'|\leq\ep_{0}
(G-\ep_{0})\sqrt{\psi} \left\{A-2\ep_{0}+
\frac{3}{4}\delta \frac{[\psi']^{2}}{\sqrt{\psi}} \right\}.$$

To this end, we fix once for all two real numbers $r_{7}$ and $r_{8}$
such that 
$$r_{5}<r_{7}<\frac{r_{5}+r_{6}}{2}<r_{8}<r_{6}.$$

When $y\in(r_{5},r_{7}]\cup[r_{8},r_{6})$ we have that $|\psi'(y)|$ is
bounded from below by a positive constant.  Therefore
$$\delta^{3/2}|k\psi'(y)|\leq\ep_{0}(G-\ep_{0})
\frac{3}{4}\delta[\psi'(y)]^{2}$$
provided that $\delta$ is small enough. When $y\in[r_{7},r_{8}]$ we
have that $\psi(y)$ is bounded from below by a positive constant, hence
in this case
$$\delta^{3/2}|k\psi'(y)|\leq
\ep_{0}(G-\ep_{0})(A-2\ep_{0})\sqrt{\psi(y)}$$
provided that $\delta$ is small enough.  This completes the proof of
(\ref{FBP2:eq-w-true}) and shows that $w$ satisfies (w7) whenever
$\delta$ is small enough.

\paragraph{\textmd{\emph{Proof of (\ref{FBP2:v>w})}}}

The argument is analogous to the proof of the corresponding
inequality in Theorem~\ref{thm:FBP1}.  Let us choose a positive value
of $\delta$ in such a way that $w$ satisfies (w1) through (w7).  Let
us assume that (\ref{FBP2:v>w}) is not true, and let us set
$$t_{0}:=\inf\left\{t\in[0,t_{\star}]:\exists r\in[r_{5}+kt,r_{6}+kt]
\mbox{ such that }v(r,t)-w(r,t)<0\right\}.$$

From (w5) we have that $t_{0}>0$.  Moreover, from the definition of
$t_{0}$ we have that
\begin{equation}
	v(r,t)-w(r,t)\geq 0
	\quad\quad
	\forall(r,t)\in\mathcal{D}_{\star}\mbox{ with }t\leq t_{0}.
	\label{FBP2:v-w-0}
\end{equation}

Finally, due to the continuity of $v$ and $w$, we deduce also that
there exists $r_{0}\in[r_{5}+kt_{0},r_{6}+kt_{0}]$ such that
$v(r_{0},t_{0})-w(r_{0},t_{0})=0$.  

We claim that $r_{0}\neq r_{5}+kt_{0}$.  Indeed let us assume by
contradiction that $r_{0}=r_{5}+kt_{0}$.  Then $w(r_{0},t_{0})=0$,
hence also $v(r_{0},t_{0})=0$.  By (v2) it follows that $r_{0}$ is a
minimum point for the function $r\to v(r,t_{0})$.  By (v4) we have
therefore that $v_{r}(r_{0},t_{0})=0$.  Keeping (w3) into account, we
deduce that
\begin{equation}
	(v-w)_{r}(r_{0},t_{0})=v_{r}(r_{0},t_{0})-w_{r}(r_{0},t_{0})<0.
	\label{FBP2:v-w-x}
\end{equation}

On the other hand, from (\ref{FBP2:v-w-0}) we know also that
$r_{0}=r_{5}+kt_{0}$ is a minimum point for the function $r\to
v(r,t_{0})-w(r,t_{0})$ restricted to the interval
$[r_{5}+kt_{0},r_{6}+kt_{0}]$.  Since the minimum point is the
left-hand endpoint of the interval, we deduce that
$(v-w)_{r}(r_{0},t_{0})\geq 0$, which contradicts (\ref{FBP2:v-w-x}).

In a symmetric way we prove that $r_{0}\neq r_{6}+kt_{0}$.  So we are
left with the case in which $t_{0}>0$ and
$r_{0}\in(r_{5}+kt_{0},r_{6}+kt_{0})$.  In this case
$0<v(r_{0},t_{0})=w(r_{0},t_{0})<c_{2}$, hence both $v$ and $w$ are
smooth in a neighborhood of this point and fulfil (v5) and (w7),
respectively.  In particular, since $r_{0}$ is always a minimum point
of the function $r\to v(r,t_{0})-w(r,t_{0})$, and now $r_{0}$ is in
the interior of the interval $(r_{5}+kt_{0},r_{6}+kt_{0})$, we have
that
\begin{equation}
	v_{r}(r_{0},t_{0})= w_{r}(r_{0},t_{0})\quad\mbox{ and }\quad
	v_{rr}(r_{0},t_{0})\geq w_{rr}(r_{0},t_{0}).
	\label{FBP2:vxx-wxx}
\end{equation}

Let us consider now the time derivatives.  On the one hand, from
(\ref{FBP2:v-w-0}) we deduce that $(v-w)_{t}(r_{0},t_{0})\leq 0$.  On the
other hand, from (v5), (w7), and (\ref{FBP2:vxx-wxx}), in the point
$(r_{0},t_{0})$ we have that
\begin{eqnarray*}
	(v-w)_{t} & > &
	g(v)\left\{v_{rr}+f(r,t,v,v_{r})+A\right\}-
	g(w)\left\{w_{rr}+f(r,t,w,w_{r})+A\right\} \\
	 & = & g(w)\left\{v_{rr}+f(r,t,w,w_{r})+A\right\}-
	g(w)\left\{w_{rr}+f(r,t,w,w_{r})+A\right\} \\
	 & = & g(w)\left( v_{rr}-w_{rr}\right)\ \geq\ 0. %\\
%	 & \geq & 0.
\end{eqnarray*}

This rules out the last case and completes the proof of
(\ref{FBP2:v>w}).

\label{NumeroPagine}

\end{document}